\documentclass{article2}
\usepackage{verbatim}
\usepackage[latin1]{inputenc}
\usepackage[T1]{fontenc}
\usepackage{amsmath,amssymb}
\usepackage{makeidx}
\usepackage[dvips]{graphicx}
\usepackage{epsfig}
\usepackage{float}
\usepackage{algorithm}
\usepackage{algorithmic}
\newtheorem{prop}{Proposition}[section]
\newtheorem{teo}{Theorem}[section]
\date{}

\author{Antonio Bernini\thanks{Dipartimento Sistemi e Informatica,
viale G.B. Morgagni 65, 50134 Firenze, Italy
\tt{bernini@dsi.unifi.it; fanti@dsi.unifi.it; ely@dsi.unifi.it}}
\and Irene Fanti$^*$ \and Elisabetta Grazzini$^*$}

\title{An exhaustive generation algorithm for Catalan
objects and others}
\begin{document}
\hyphenation{generated} \maketitle
\begin{abstract}
In this paper we present a CAT generation algorithm for Dyck paths
with a fixed length $n$. It is the formalization of a method for
the exhaustive generation of this kind of paths which can be
described by means of two equivalent strategies. The former is
described by a rooted tree, the latter lists the paths by means of
three operations which, as we are going to see, are equivalent to
visit the tree. These constructions are strictly connected with
ECO method and can be encoded by a rule, very similar to the
\emph{succession rule} in ECO, with a finite number of labels for
each $n$. Moreover with a slight variation this method can be
generalized to other combinatorial classes like Grand Dyck or
Motzkin paths.
\end{abstract}

\section{Introduction}
One of the most important aims in combinatorics has always been
the generation of objects of a particular class according to a
fixed parameter. Actually many practical questions require, for
their solution, an exhaustive search through out all the objects
in the class. In general, the idea is to find methods to list in a
particular order combinatorial objects without either repetitions
or omissions so that it's possible to deduce a recursive
construction of studied class. We are talking about exhaustive
generation algorithms \cite{1} which can be seen as an enumerating
techniques where each object is counted and recorded once as it is
generated \cite{3}. Often, these algorithms are useful in diverse
areas such as hardware and software testing, thermodynamic,
biology and biochemistry (\cite{4,2,5}) where it could be helpful
to have a particular order of the objects.

In literature one of the common approach has been the generation
of the combinatorial elements in such a way that two successive
items differ only slightly; in this sense a well known example is
the classical binary reflected Gray code scheme for listing
$n$-bit binary numbers so that successive items differ in exactly
one bit position. Gray codes find a lot of applications in many
different areas (for more details and examples see
\cite{11,12,20,6,9,7,13,15,19,18,17,14,8,10,16}).

In this paper we present a method which allows us to generate all
and only objects of a combinatorial class, fixed the size. We
focus attention on Dyck paths and we introduce two strategies for
generating them. The former uses an operator which can be
described by a rooted tree, the latter lists the objects using
three operations and it corresponds to visit all nodes of the
tree. Both of them, as we will see, use only a constant amount of
computations per object in amortized sense and so they have CAT
property \cite{22}. These methods are efficient, in fact it is
well known that the primary performance goal in an algorithm for
listing a combinatorial family is to have a running time
proportional to the number of objects produced \cite{26}.

Moreover our method is similar to ECO method and the analogy
suggests to look for a rule, similar to the \emph{succession rule}
in ECO, for encoding the construction of the list (for more
details about ECO method and \emph{succession rule} see \cite{22,
21, 27, 23, 24}). Nevertheless the previous approaches (see in
\cite{22}) used string of integers for encoding the objects and it
was necessary to have an other algorithm, which required a
proportional amount of computations, to come back to the objects.
On the contrary our method directly uses the combinatorial objects
and it generates the paths by means of operations with a constant
cost. In section \ref{def} we give some preliminaries and
notations 
while, in the other sections, we present our main idea.

\section{Preliminaries and notations}\label{def}In this section we
give some notations which are necessary to introduce our
method.

\begin{figure}[!ht]
\begin{centering}
\includegraphics[width=0.6
\textwidth]{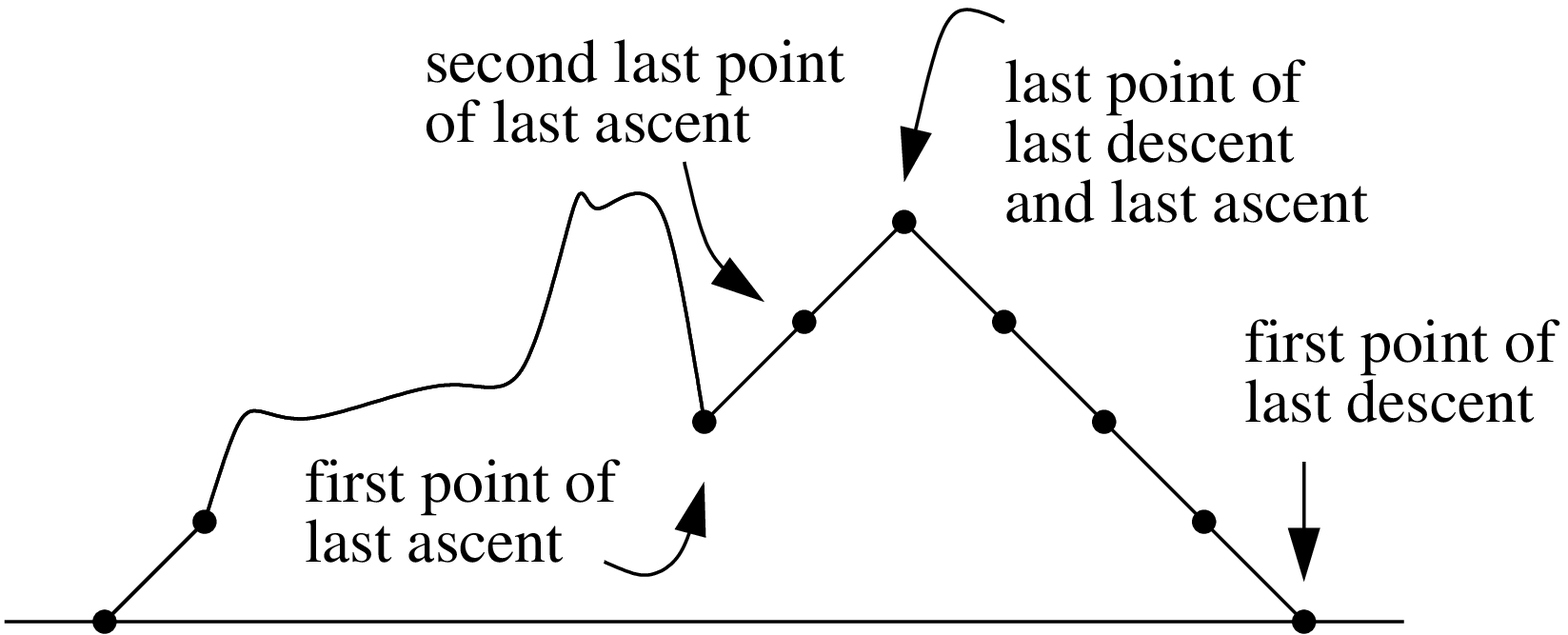} \caption{Numeration of points of Dyck path's
last descendent and ascent.} \label{tr:num}
\end{centering}
\end{figure}

As we said above, we consider the combinatorial class of Dyck
paths and we give some definitions useful for our work. We define
a \emph{path} like a sequence of points in $\mathbb{N} \times
\mathbb{N}$ (they have never negative coordinates) and a
\emph{step} like a pair of two consecutive points in the path. A
\emph{Dyck path} is a path $\mathcal{D} : = \{ s_0, s_1, \ldots,
s_{2n} \}$ such that $s_0 = ( 0 , 0 )$ , $s_{2n} = ( 2n, 0 )$ and
only having northeast ($s_i = ( x , y )$, $s_{i + 1 } = ( x + 1 ,
y + 1 )$) or southeast ($s_i = ( x , y )$, $s_{i + 1 } = ( x + 1 ,
y - 1 )$) steps; so the number of northeast steps is equal to the
number of southeast steps and we can define \emph{path's length}
the number of its steps. In particular $\mathcal{D}_n$ is the set
of Dyck paths with length $2n$ i.e with $n$ northeast steps. In
the sequel we say that $\mathcal{D}\in\mathcal{D}_n$ has
\emph{size} $n$.
\newline A \emph{peak} (resp.\emph{valley}) is a point $s_i$ such
that step ($s_{i - 1}$, $s_i$) is a northeast (southeast) and the
step ($s_i$, $s_{i + 1 }$) is a southeast (northeast); moreover we
say \emph{pyramid} $p_h \ , \forall \ h \in \mathbb{N}$, a
sequence of $h$ northeast steps following by $h$ southeast steps
such that if $(s_i, s_{i + 1 })$ is the first northeast step and
$(s_{i+2h-1}, s_{i+2h})$ is the last southeast of this sequence,
then $s_i = (x, 0)$ and $s_{i+2h} = (x+2h, 0)$. We also define
\emph{last descent} (\emph{ascent}) the southeast (northeast)
steps' last sequence of a Dyck path and we conventionally number
its points from right (left) to left (right); clearly the last
point of last descent always coincides with last point of last
ascent (see Figure \ref{tr:num}).
Moreover, if we say \emph{height} $h(s_i)$ of a point $s_i$ its
ordinate and \emph{non-decreasing point} the extremity $s_{i + 1}$
of a northeast step ($s_i$, $s_{i + 1 }$), then we can define
\emph{area of a path} the sum of its non-decreasing points'
heights and \emph{maxima area path} $P ^n_{max}$ the pyramid $p_n$
that contains, in geometric sense, all the paths of its size.
Finally we call a path $P$ ``active'' if we obtain another Dyck
path when the first and the last step of $P$ are taken off. This
is equivalent to say that $P$ does not have valleys with height
$h=0$.

Given a class of combinatorial objects $\mathcal{C}$ and a
parameter $\gamma : \mathcal{C}$ $\longrightarrow$ $\textbf{N}^+$
such that ${\mathcal{C}}_n = \{x\in {\mathcal{C}} : \gamma \ ( x )
= n\}$ is a finite set for all $n$, we define a \emph{generating
tree}. We assume there is only one element of minimal size in
$\mathcal{C}$ and we describe the recursive construction of this
set by using a rooted tree in which each node corresponds to an
object. In particular, the vertices on the $n$th level represent
the elements of $\mathcal{C}_n$, the root of the tree is the
smallest element and the branch, leading to the node, encodes the
choices made in the construction of the object. Starting from this
idea and choosing the combinatorial class $\mathcal{D}$ of Dyck
paths, we introduce another kind of generating tree which
describe, fixed the size $n$, the recursive construction of
$\mathcal{D}_n$. In the sequel we denote it with
$\mathcal{D}_n$\emph{-tree} which clearly has a finite number of
levels and each object has the same size, regardless of the level.
The structure of a generating tree can be elegantly describe by
means of the notion of \emph{succession rule}. Moreover, as we
have just said, our algorithm is based on the ECO method which is
a general method to enumerate combinatorial objects. The basic
idea of this one is the definition of a recursive construction for
$\mathcal{C}$ by means of an operator $\phi$ which performs a
``local expansion'' on the objects (for more details see
\cite{21}).
\section{Dyck paths}\label{met}

\begin{figure}[!ht]
\begin{centering}
\includegraphics[width=0.9\textwidth]{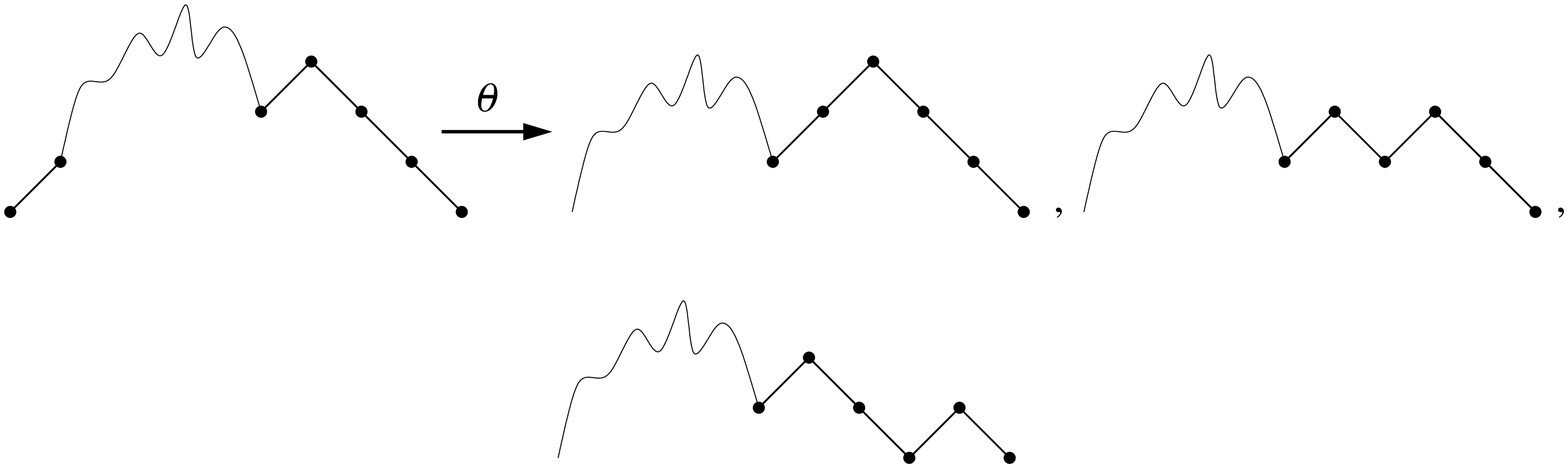}\caption{} \label{theta}
\end{centering}
\end{figure}

\begin{figure}[!ht]
\begin{centering}
\includegraphics[width=0.9\textwidth]{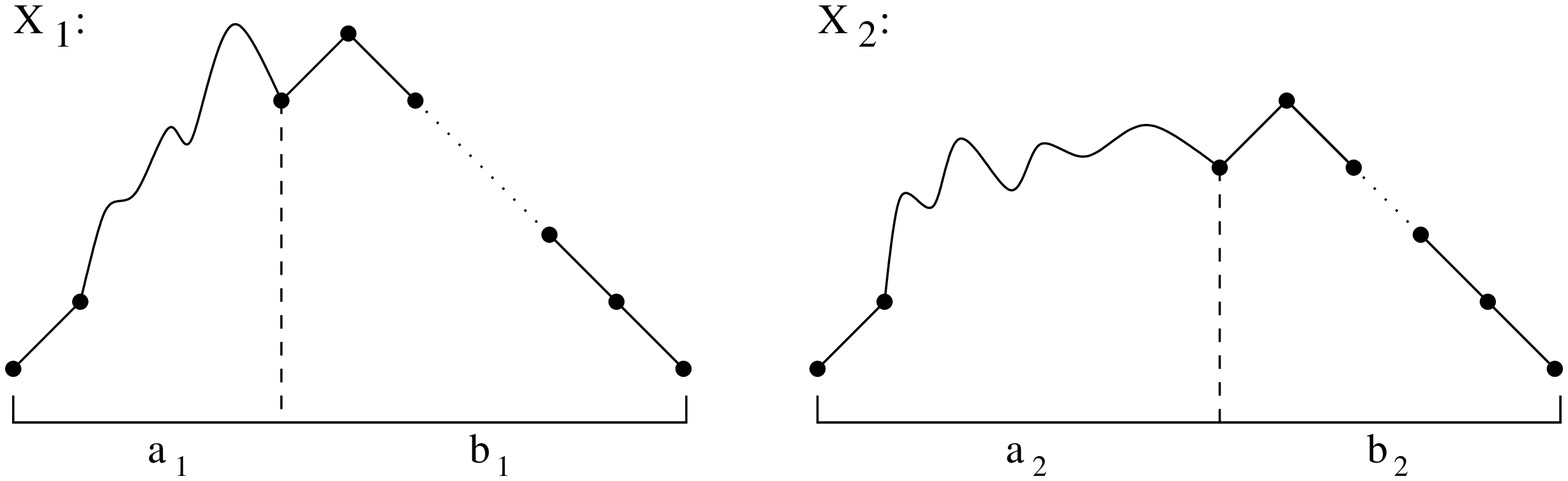}\caption{}
\label{tr:x}
\end{centering}
\end{figure}

We start to define an operator which constructs $\mathcal{D}_n$;
we study this operator for $n\geqslant3$ being cases $n=1$ and
$n=2$ trivial.\vspace{0.5cm}\newline\textbf{$\theta$ OPERATOR:}
\begin{enumerate}
    \item Consider $P^n_{max}$ like the first path.
    \item Take off the first and the last path's step and
    insert a peak in every point of the obtained path's last descent except for the last point.
    Every insertion generates a new Dyck path.
    \item For each new generated path repeat the following actions until active paths are generated:
    \begin{enumerate}
        \item [3.1] take off the first and the last path's step
        \item [3.2] insert a peak in every point of the obtained path's last descent.
                    Every insertion generates a new path.
    \end{enumerate}
\end{enumerate}
In Figure \ref{theta} we give an example of $\theta$ operator's
action.
\newline We prove that $\theta$ satisfies the following conditions:
\begin{prop}\label{p1}\mbox{}
\begin{itemize}
    \item [1.] $\forall$ $X_1$ , $X_2 \in \theta(P^n_{max})$, then $X_1\neq
    X_2$;
    \item [2.] $\forall$ $X_1$ , $X_2 \in \mathcal{D}_n$ and $X_1\neq X_2$, then
    $\theta (X_1) \cap \theta(X_2) = \emptyset$.
\end{itemize}
\end{prop}
\begin{prop}\label{p2}$\forall$ $Y \in \mathcal{D}_n$ $\exists$ a finite succession $X_0, X_1, \ldots, X_k$
with  $k \in \mathbb{N}$ and $X_k=Y$ such that :\begin{itemize}
    \item $X_0=P^n_{max}$;
    \item $X_{i+1} \in \theta(X_i)$ \hspace{0.5cm}$0\leq$ i $\leq$
    $k-1$.
\end{itemize}
\end{prop}
\textbf{Proof Proposition \ref{p1}}: We prove point 2 since point
1 of the proposition is trivial. Consider $X_1$ , $X_2 \in
\mathcal{D}_n$, $X_1\neq
    X_2$ 
    and divide both  $X_1$ and $X_2$ in two parts as shown in Figure
    \ref{tr:x}.
   \newline If $b_1 \neq b_2$, they remain distinct after the application
   of $\theta$ since it operates just on these parts. On the other
   hand if $b_1 = b_2$, then $a_1\neq a_2$ and after the
   application of $\theta$ $a_1$ and $a_2$ remain different, then $\theta (X_1) \cap \theta (X_2) =
   \emptyset$ in both cases.
\vspace{0.3cm}\newline \textbf{Proof Proposition \ref{p2}}: We
consider a general path Y and we apply the inverse of $\theta$
operator on it; clearly $\theta^{-1}$ operator takes off the
righter peak of $Y$ and inserts a northeast step at the beginning
of the path and a southeast step at the end. We have two
possibilities:
\begin{description}
    \item[1.] The last ascent of $Y$ has only one step, so in the
    obtained path $\theta^{-1}(Y)$ the peaks' number is reduced by
    one.
    \item[2.] The last ascent of $Y$ has at least two steps, so
    the number of last ascent's steps in $\theta^{-1}(Y)$ is
    reduced by one.
\end{description}
It is clear that after $k$ times, for $k\in \mathbb{N}$, the
number of peaks in $\theta^{-k}(Y)$ is one and
$\theta^{-k}(Y)=P^n_{max}$.
\begin{figure}[!ht]
\begin{centering}
\includegraphics[width=0.95\textwidth]{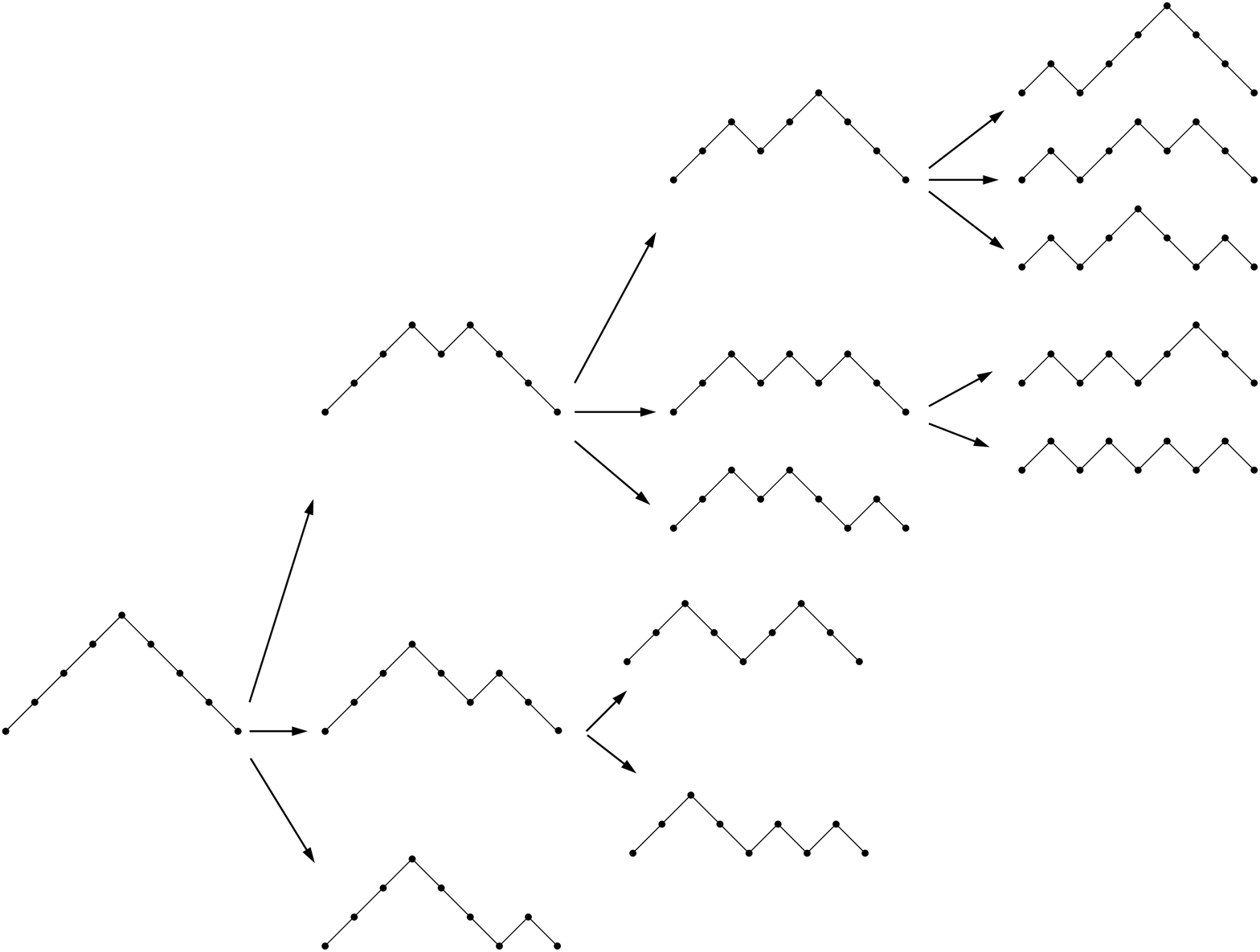} \caption{$\mathcal{D}_4$-tree.}
\label{tr:alb1}
\end{centering}
\end{figure}
\newline Now we pass to describe $\theta$'s construction by
using a rooted
tree:\vspace{0.5cm}\newline\textbf{$\mathcal{D}_n$-tree ROOTED
TREE:}
\begin{enumerate}
    \item The root is $P^n_{max}$ and it is at level zero;
    \item if $X \in \mathcal{D}_n$-tree is at level $k\geq0$ then $Y \in \theta
    (X)$ is a son of $X$ and it is at level $k+1$.
\end{enumerate}
In Figure \ref{tr:alb1} $\mathcal{D}_4$-tree is illustrated.
\begin{teo}$\mathcal{D}_n=\mathcal{D}_n$-tree
\end{teo}
\textbf{Proof:} Given $X\in\mathcal{D}_n$-tree; it is clear that
$X$ is a Dyck path. Moreover, Proposition \ref{p1} assures there
are not two copies of the same path in $\mathcal{D}_n$-tree
$\Rightarrow |\mathcal{D}_n$-tree$|\leq |\mathcal{D}_n|$ \
$\Rightarrow \mathcal{D}_n$-tree $\subseteq \mathcal{D}_n$.
\newline Vice versa given $Y \in \mathcal{D}_n$, Proposition \ref{p2} assures that it
is always possible to find a finite succession which joins
$P^n_{max}$ path to $Y$; so $Y \in \mathcal{D}_n$-tree since
$P^n_{max}$ is in $\mathcal{D}_n$-tree \ $\Rightarrow
\mathcal{D}_n \subseteq \mathcal{D}_n$-tree.
\subsection{Succession rule}

\begin{figure}[!ht]
\begin{centering}
\includegraphics[width=1.0
\textwidth]{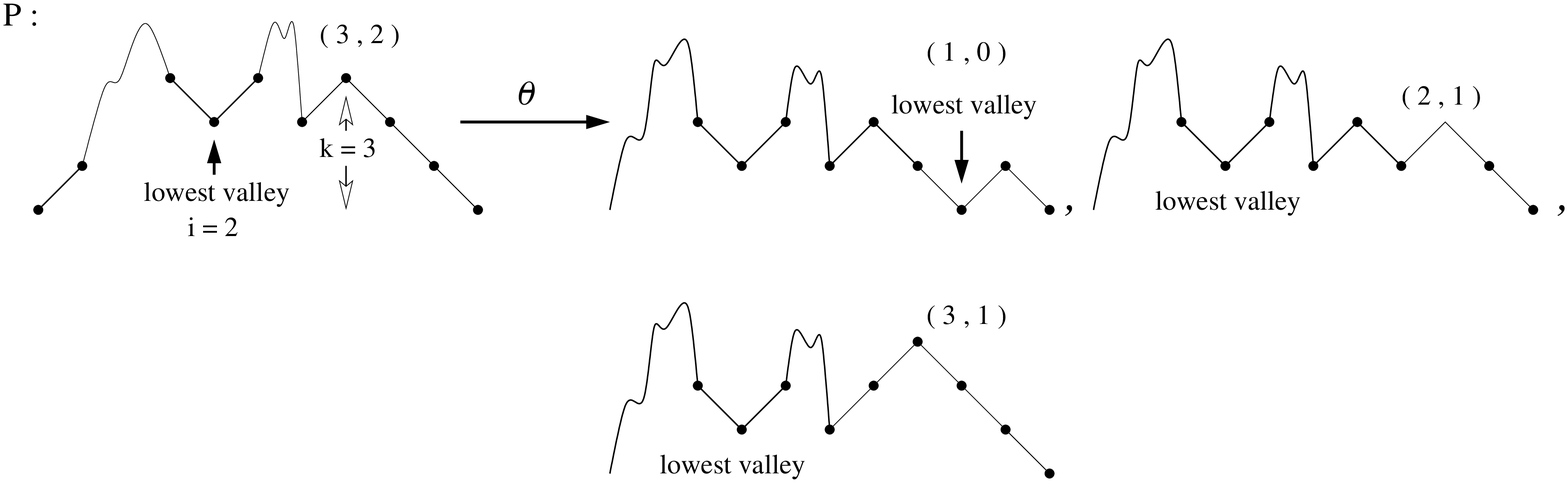}\caption{} \label{theta4}
\end{centering}
\end{figure}

Now we give a succession rule to describe $\mathcal{D}_n$-tree. As
we have said above, given a path $P$, we have $\theta(P)\neq
\emptyset$ if and only if it is active, i.e. if it has not valleys
with height $h=0$. Moreover, from the definition of $\theta$
operator it is clear that the number of a path's sons is equal to
the number of steps in its last descent. So, we have to label each
path with an information which says us the number of its sons and
the height of its lowest valley. We use the following notation to
connect the label of a parent $P$, having the height of its lowest
valley equal to $i$, with the labels of its $k$ sons:
    $$
    (k,i)\hookrightarrow(c_1)(c_2) \ldots (c_k).
    $$
Moreover, each of these $k$ paths has the last descent with length
$s$, with $s = 1, 2, \ldots , k$. Now $\theta$ operator, after
having taken off the first and the last step of $P$, inserts a
peak in one of the last descent's point of the obtained path
$\overline{P}$. This insertion increases the number of valleys in
the generated path by one, with the exception of $Y$ obtained by
inserting the peak in the last point of $\overline{P}$'s last
descent, since in this case, the path has the same number of
valleys of its father $P$. So, it is clear that the height $j$ of
generated paths' lowest valley depends of the insertion of the
peak. Indeed, if $\theta$ inserts the peak in the $t$-th point of
the $\overline{P}$'s last descent with $1\leqslant t \leqslant
i-1$, then $j=t-1$, i.e. the lowest valley is generated by the
peak insertion. On the other hand if $i\leqslant t \leqslant k$,
then $j=i-1$, i.e. the lowest valley is the same of
$\overline{P}$. In Figure \ref{theta4} we give an example of
$\theta$'s action on a path with label $(3,2)$.

\begin{figure}[!ht]
\begin{centering}
\includegraphics[width=1.0
\textwidth]{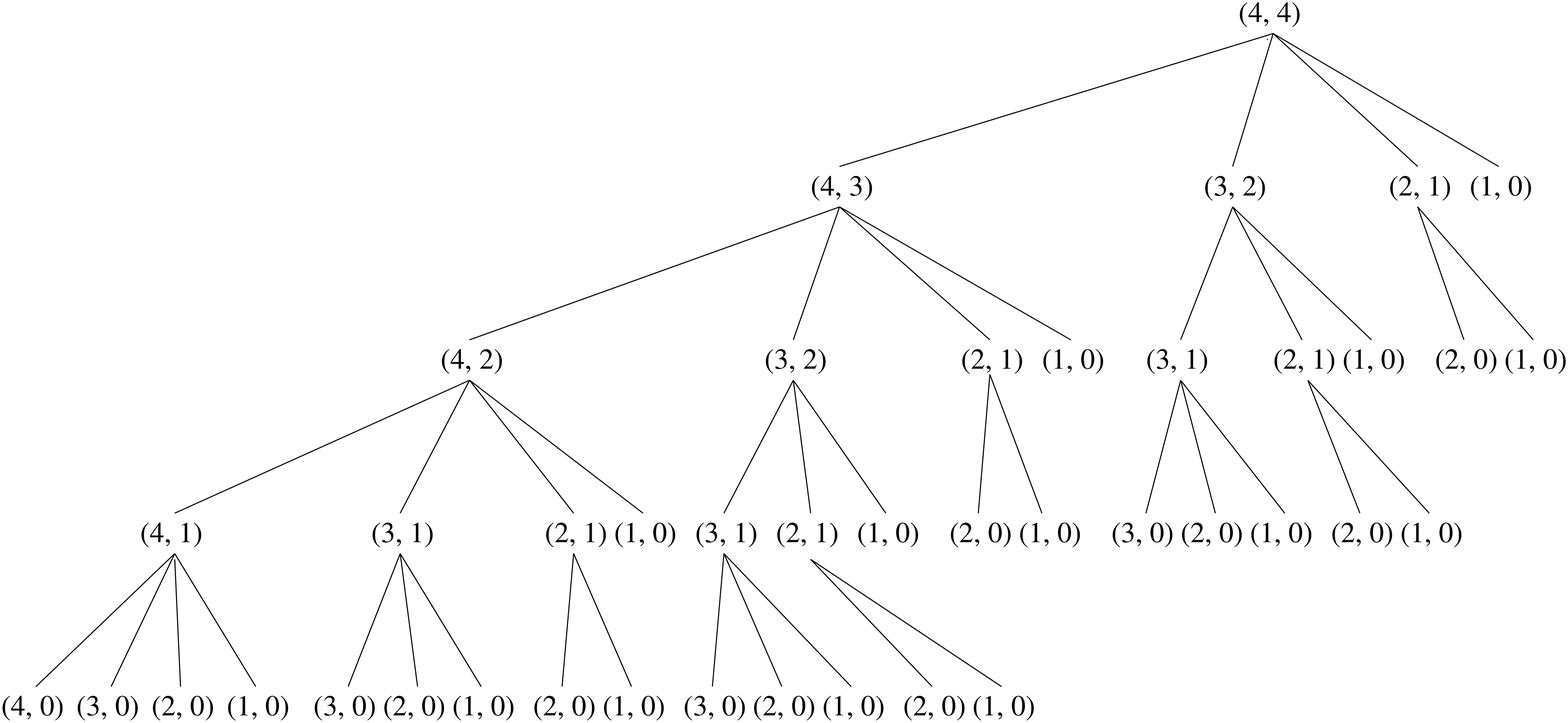}\caption{} \label{tr:eti}
\end{centering}
\end{figure}

We can give the production:
$$
(k,i)\hookrightarrow(1,0)(2,1) \ldots (i,i-1)(i+1,i-1) \ldots (k,
i-1).
$$
We notice that the root of $\mathcal{D}_n$-tree doesn't have
valleys and the second index of its label could be empty;
nevertheless, for convenience, we have decided to label the root
by $(n-1,n-1)$. Finally, we obtain the following succession rule:
$$\left \{
\begin{array}{ll}
(n-1, n-1)\\
(k,i)\hookrightarrow(1,0)(2,1) \ldots (i,i-1)(i+1,i-1) \ldots (k,
i-1)\\
\end{array} \right.$$
It is clear that labels with $i=0$ correspond to paths with at
least a valley with height $h=0$ and they do not generate any
other path by $\theta$ operator. In Figure \ref{tr:eti} we give an
example of generating tree for $n=5$.

\subsection{The generating algorithm} 
\begin{figure}[!ht]
\begin{centering}
\includegraphics[width=1.02
\textwidth]{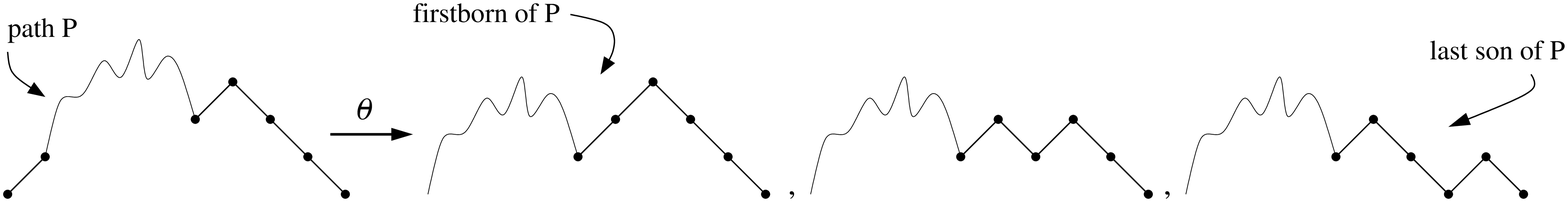}\caption{} \label{tr:es}
\end{centering}
\end{figure}

\begin{figure}[!ht]
\begin{centering}
\includegraphics[width=0.85
\textwidth]{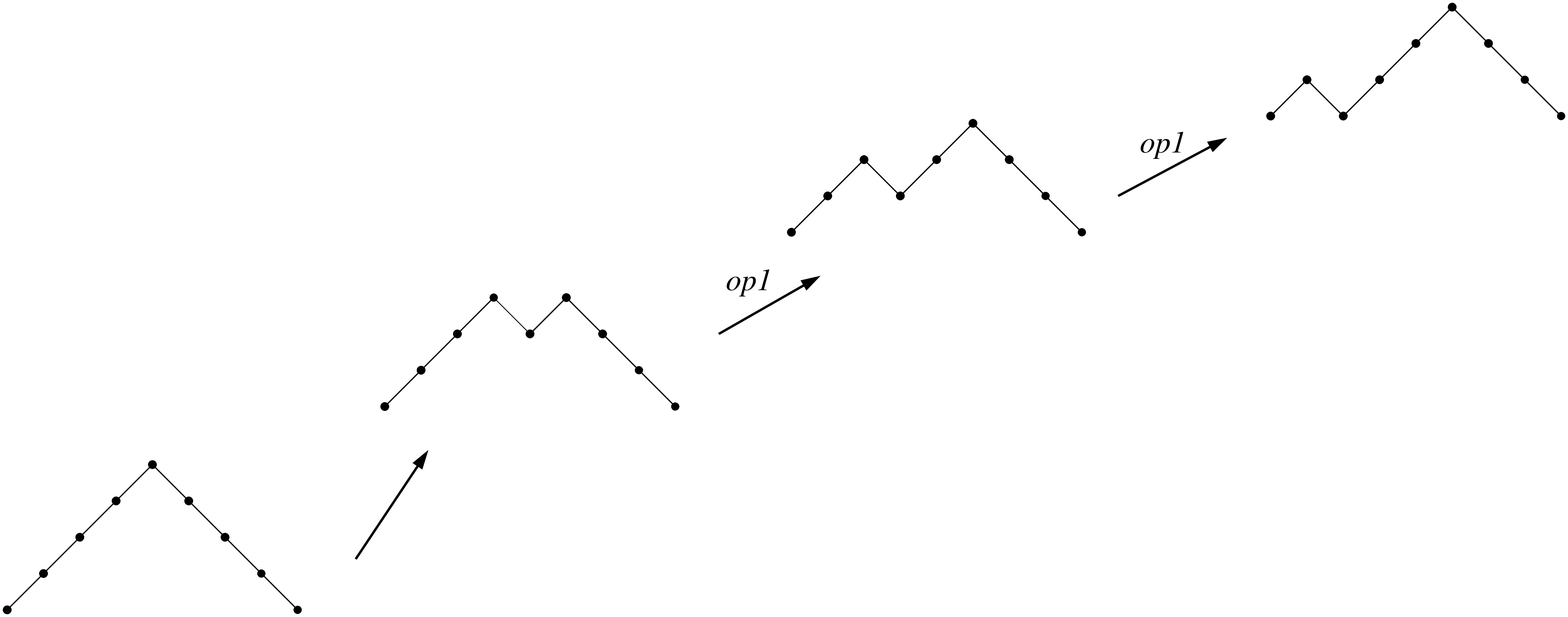}\caption{} \label{tr:albop1}
\end{centering}
\end{figure}

\begin{figure}[!ht]
\begin{centering}
\includegraphics[width=0.85\textwidth]{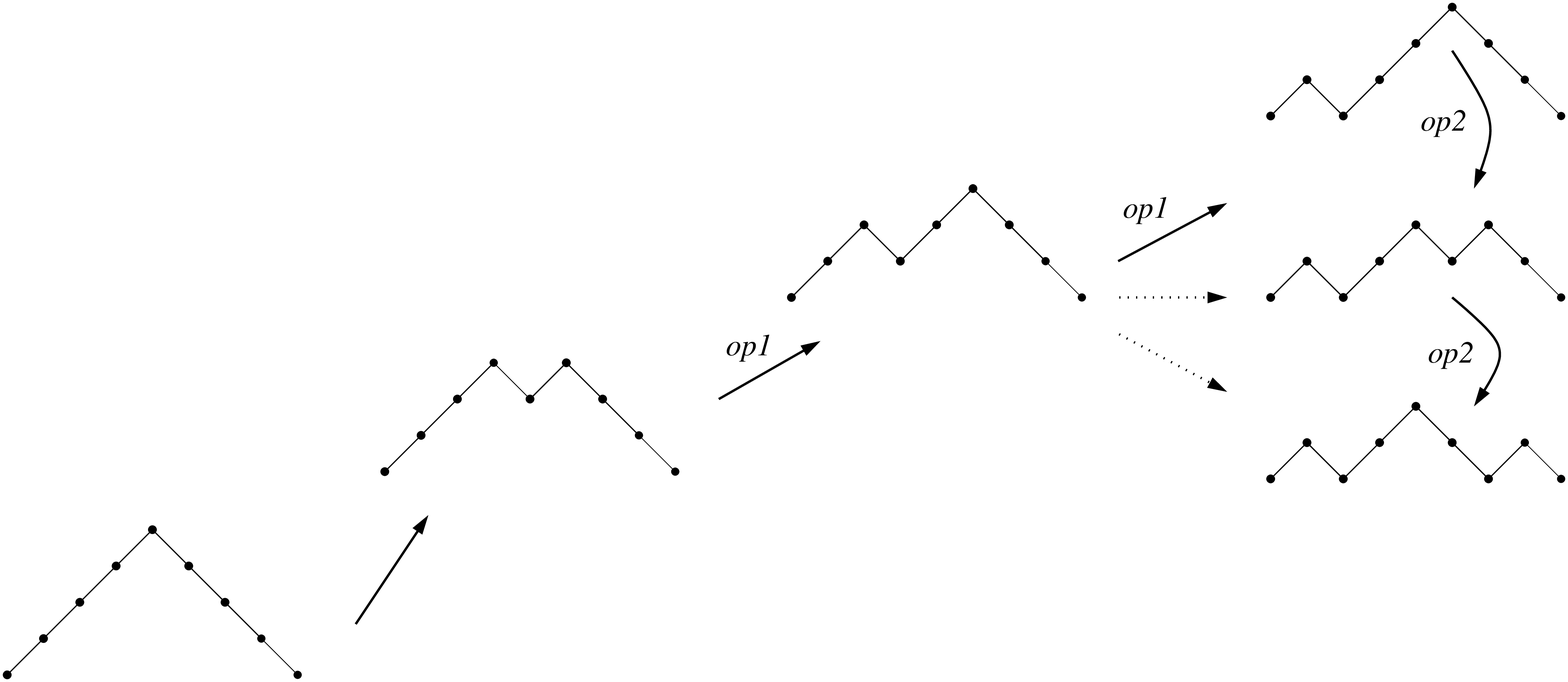}\caption{}
\label{tr:albop2}
\end{centering}
\end{figure}

As we have just seen, $\theta$ operator can be
described by a rooted tree and $\mathcal{D}_n$'s paths are
generated according to the $\mathcal{D}_n$-tree's levels.
Nevertheless we wish to find a method which sequentially lists the
objects so that everyone is generated only by the last generated
path. This operation corresponds to visit all the nodes of
$\mathcal{D}_n$-tree and for this reason it's helpful to order the
sons of $X$ path according to the decreasing length of their last
descent so that the last one ends in $p_1$. In particular, the
last $P^n_{max}$'s son is made by $p_{n - 1 }$ followed by $p_1$.
We name ``firstborn'' of a path $P$ the son which has the longest
last descent (In Figure \ref{tr:es} we give an example of a path's
``firstborn'').
\newline Clearly the ``firstborn'' of $P^n_{max}$ can be generated simply
overturning its peak. Then we generate all ``firstborn'' paths on
the longest branch of $\mathcal{D}_n$-tree applying $(n-2)$ times
the following operation:\vspace{0.3cm}\newline \textbf{\emph{
op1:}} Take off the first and the last path's step, then insert a
peak in the last point of the last descent (see Figure
\ref{tr:albop1}).
\newline When \emph{op1} is no more applicable i.e. when we arrive at a
leaf, we proceed to generate the leaf's brothers following the
order given at the beginning of this subsection. So it is
sufficient to apply the following operation on the last generated
path: \vspace{0.3cm}\newline \textbf{\emph{op2:}} Overturn the
rightmost peak in the path (see Figure \ref{tr:albop2})
since, if $Y_i \in \theta (X)$ with $1\leqslant i \leqslant k-1$
and $k=|\theta(X)|$, then \emph{op2}$(Y_i) = Y_{i + 1 } \in \theta
(X)$. Indeed, $Y_{i+1}$ is generated from $X$ by means of $\theta$
taking off the first and the last step and inserting a peak in the
$(k - i + 2)th$ point of $X$'s last descent; the generation of
$Y_{i+1}$ can be obtained also overturning the rightmost peak in
$Y_i$.

\begin{figure}[!ht]
\begin{centering}
\includegraphics[width=0.85\textwidth]{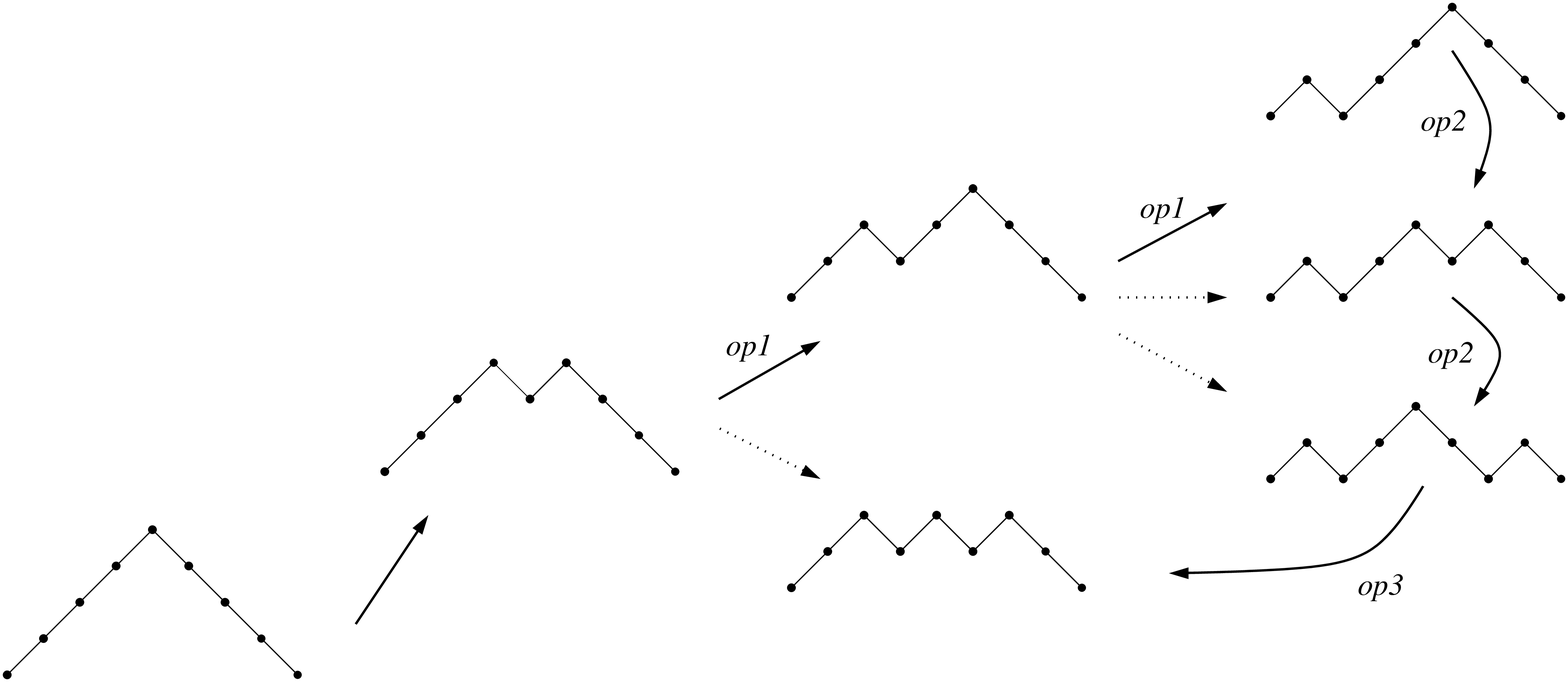}
\caption{} \label{tr:albop3}
\end{centering}
\end{figure}

\begin{figure}[!ht]
\begin{centering}
\includegraphics[width=0.9\textwidth]{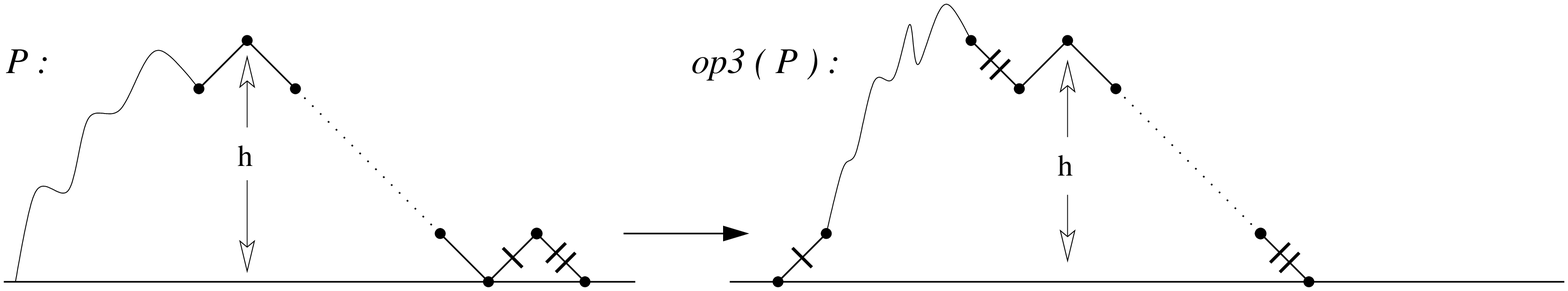}
\caption{} \label{tr:3}
\end{centering}
\end{figure}

\begin{figure}[!ht]
\begin{centering}
\includegraphics[width=0.9\textwidth]{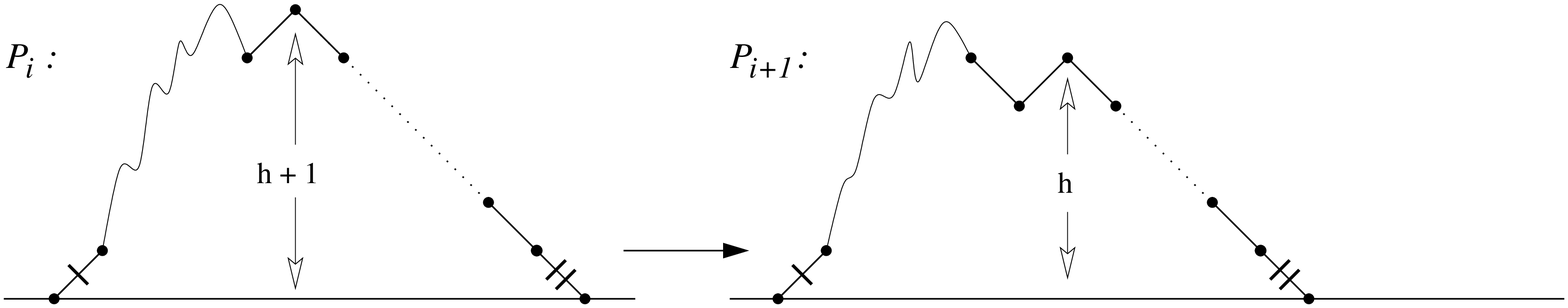}
\caption{} \label{tr:3d}
\end{centering}
\end{figure}

After the last son of $X$ is generated, we should go back to the
$\mathcal{D}_n$-tree's preceding level, in other words we should
pass to the immediately next brother of $X$, if it exists. We use
the following operation to generate the ``uncle'' of last obtained
path:\vspace{0.3cm}\newline \textbf{\emph{op3:}} Take off the
rightmost $p_1$; then insert a northeast step at the beginning of
the path and a southeast step in the second-last point of the last
ascent (see Figure \ref{tr:albop3}).
\newline \emph{Op3} allows us to pass from a path ending in
$p_1$ to its ``uncle''; this fact it's very important because we
can pass to another subtree of $\mathcal{D}_n$-tree, where we can
apply \emph{op1} and \emph{op2} again.
\newline The effects of \emph{op3} on a path $P$ ending in $p_1$
 are illustrated in Figure \ref{tr:3}.
\newline Let $P$ the last son of a path $P_i$; the
path $P_{i+1}$, $P_i$'s brother, is obtained simply overturning
its last peak (see Figure \ref{tr:3d}).
\newline As we can see, the second path in the first figure is
equivalent to $P_{i+1}$ so \emph{op3}$(P)=P_{i+1}$. \newline From
their definition \emph{op1}, \emph{op2} and \emph{op3} form a
method to visit all the nodes of $\mathcal{D}_n$-tree and so, they
generates all $\mathcal{D}_n$ paths.

\begin{figure}[!ht]
\begin{centering}
\includegraphics[width=1\textwidth]{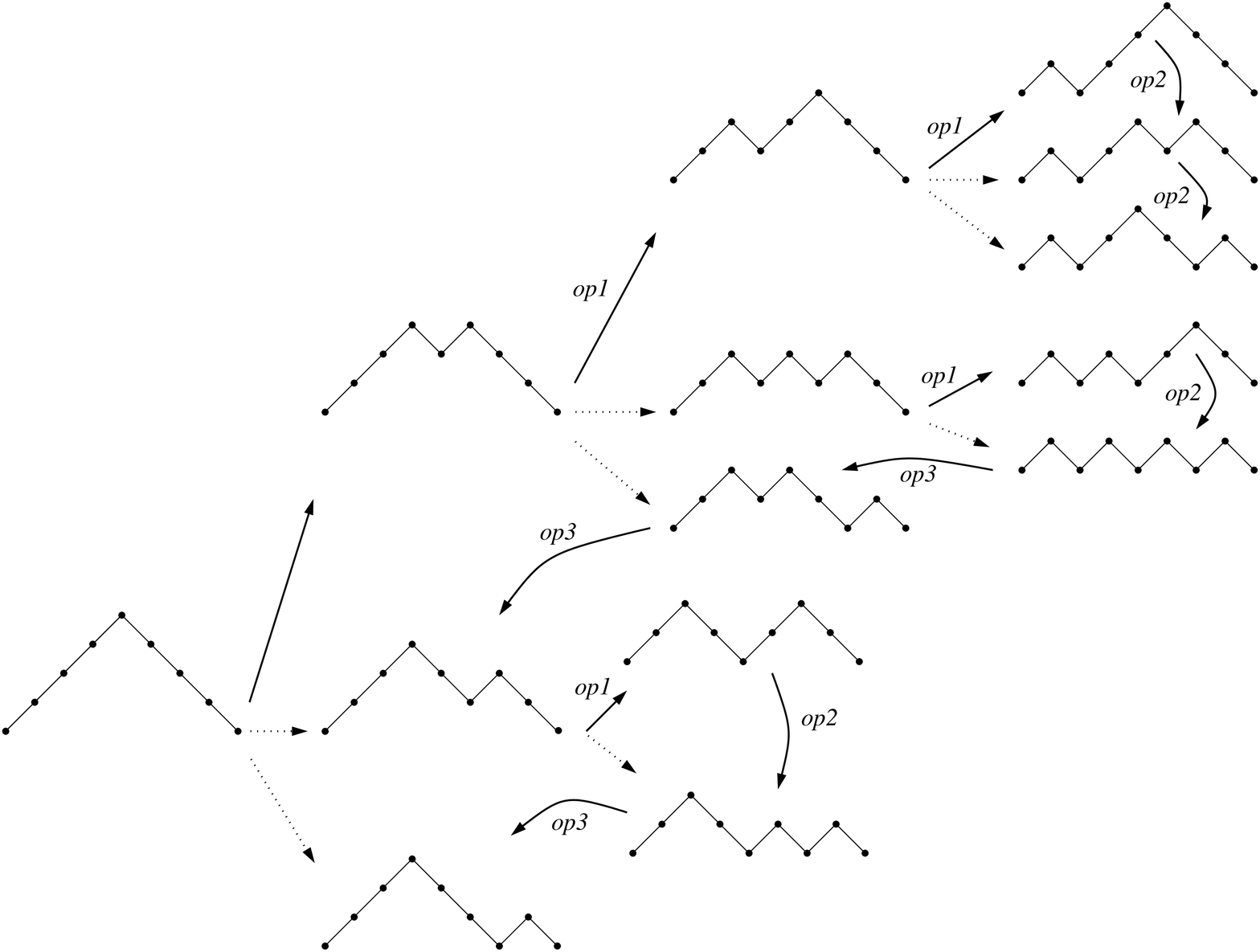}
\caption{}\label{tr:alb3op}
\end{centering}
\end{figure}

We proceed to express all that we have exemplified by means of the
following algorithm:\floatname{algorithm}{Algorithm}
\begin{algorithm}
\caption{}\label{algo}
\begin{algorithmic}
\STATE $start \hspace{0.2cm} with \hspace{0.2cm} P^n_{max}$;
 \STATE $generate \hspace{0.2cm} the \hspace{0.2cm} firstborn \hspace{0.2cm}
 son \hspace{0.2cm} of \hspace{0.2cm} P^n_{max} \hspace{0.2cm} overturning \hspace{0.2cm} its \hspace{0.2cm} peak$;
 \STATE $P:= firstborn \hspace{0.2cm}
 son \hspace{0.2cm} of \hspace{0.2cm} P^n_{max}$;
 \WHILE {$P \hspace{0.2cm} \neq \hspace{0.2cm} the \hspace{0.2cm} last
 \hspace{0.2cm} son \hspace{0.2cm} of \hspace{0.2cm} P^n_{max}$}
   \IF {it's\hspace{0.1cm}possible}
      \STATE $P^\prime:=$ \emph{op1}$(P)$
      \ELSIF {it's\hspace{0.1cm}possible}
      \STATE $P^\prime:=$ \emph{op2}$(P)$
      \ELSE
      \STATE $P^\prime:=$ \emph{op3}$(P)$
   \ENDIF;
   \STATE $P:= P^\prime$;
 \ENDWHILE
 \end{algorithmic}
\end{algorithm}
\newline\textbf{\emph{Remark:}} Observing Figure \ref{tr:alb3op}
we can notice that it's possible to have more consecutive
operations of the same kind but in particular we can have at most
two consecutive applications of \emph{op3}. Indeed we can have
only two possibilities:
\begin{description}
    \item[a)] \textbf{The path ends in $p_1$ which is preceded by a peak with height
    h~$\geq$~2} \newline \emph{Op3} works only one time because its
    application, as we can see in Figure \ref{tr:mag1},
    generates a path that has the last peak with height
    h $\geq2$.
    \begin{figure}[!ht]
\begin{centering}
\includegraphics[width=0.9\textwidth]{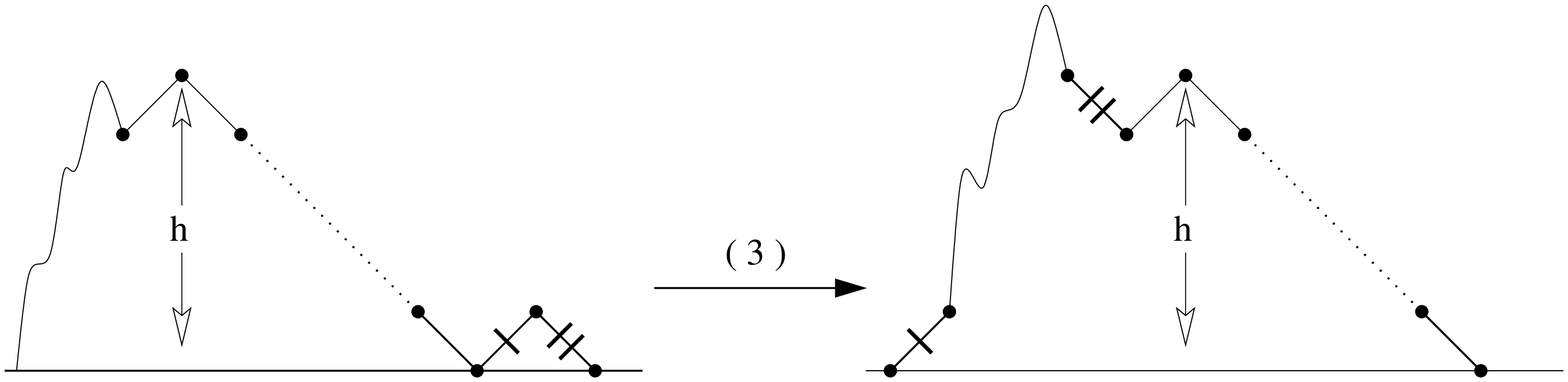}
\caption{}\label{tr:mag1}
\end{centering}
\end{figure}
    \item[b)] \textbf{The path ends in at least two $p_1$}
    \newline In this case the application of \emph{op3} generates
    a path that ends again in $p_1$; we are in case
    \textbf{a)} and the application of \emph{op3} is possible only another
    time (see Figure \ref{tr:1}).
    \begin{figure}[!ht]
\begin{centering}
\includegraphics[width=1.0\textwidth]{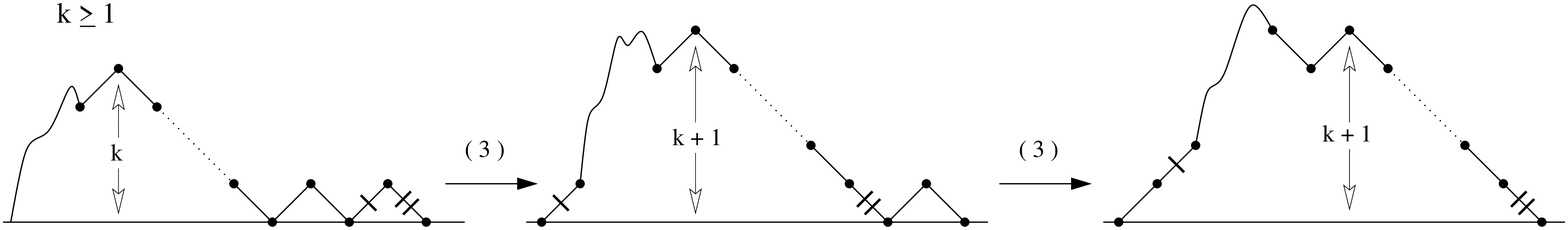}
\caption{} \label{tr:1}
\end{centering}
\end{figure}
\end{description}

\subsection{Analysis of the algorithm}
\begin{figure}[!ht]
\begin{centering}
\includegraphics[width=0.8\textwidth]{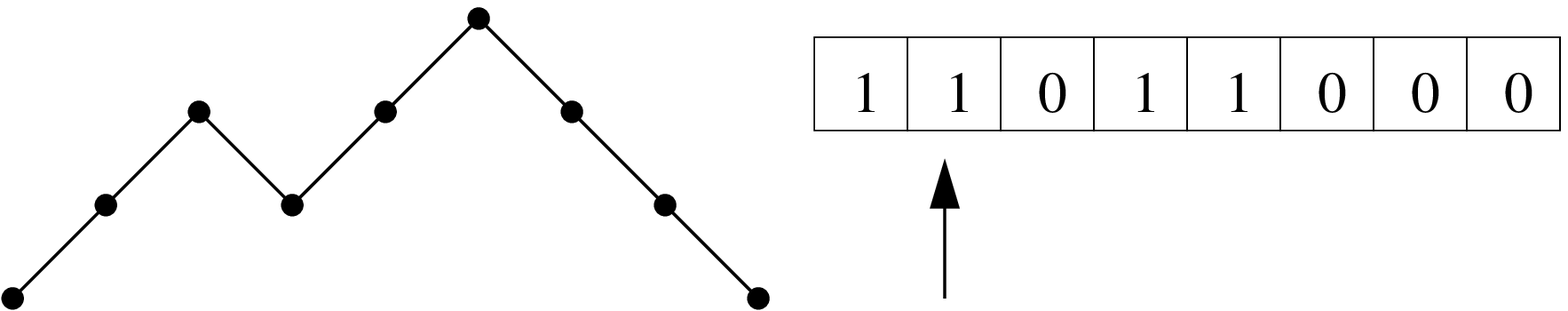}
\caption{} \label{tr:vet}
\end{centering}
\end{figure}

\begin{figure}[!ht]
\begin{centering}
\includegraphics[width=0.7\textwidth]{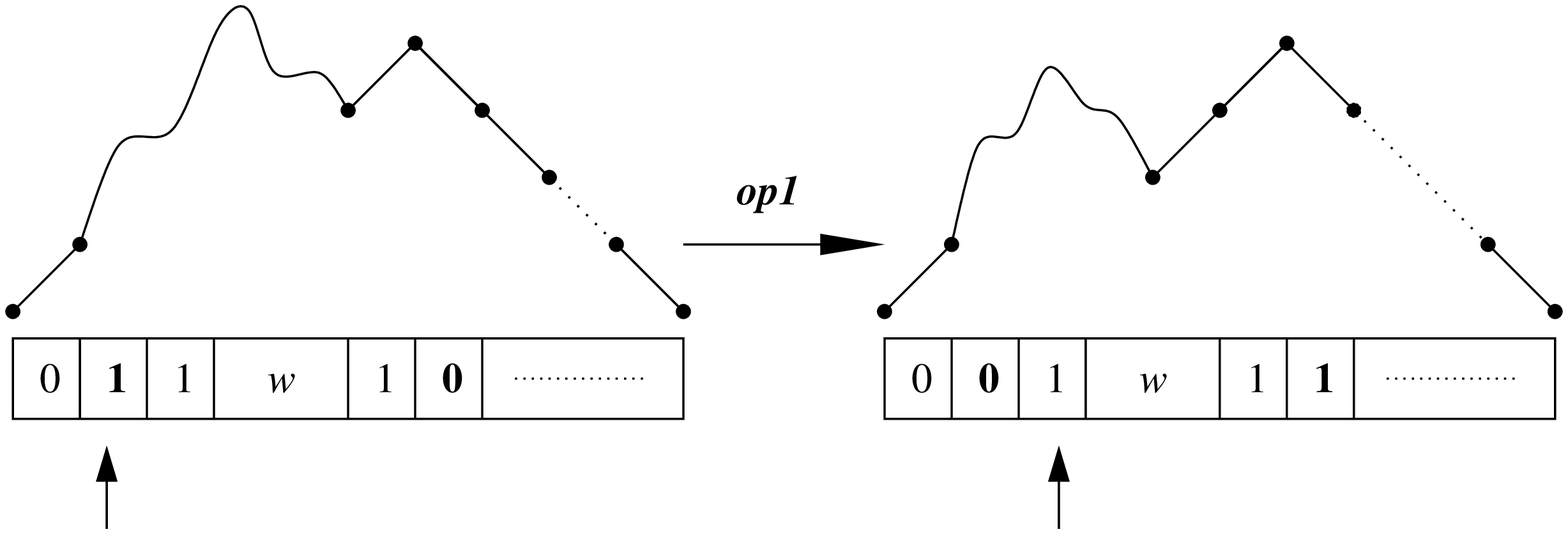}\caption{}
\label{tr:op1v}
\end{centering}
\end{figure}

\begin{figure}[!ht]
\begin{centering}
\includegraphics[width=0.7\textwidth]{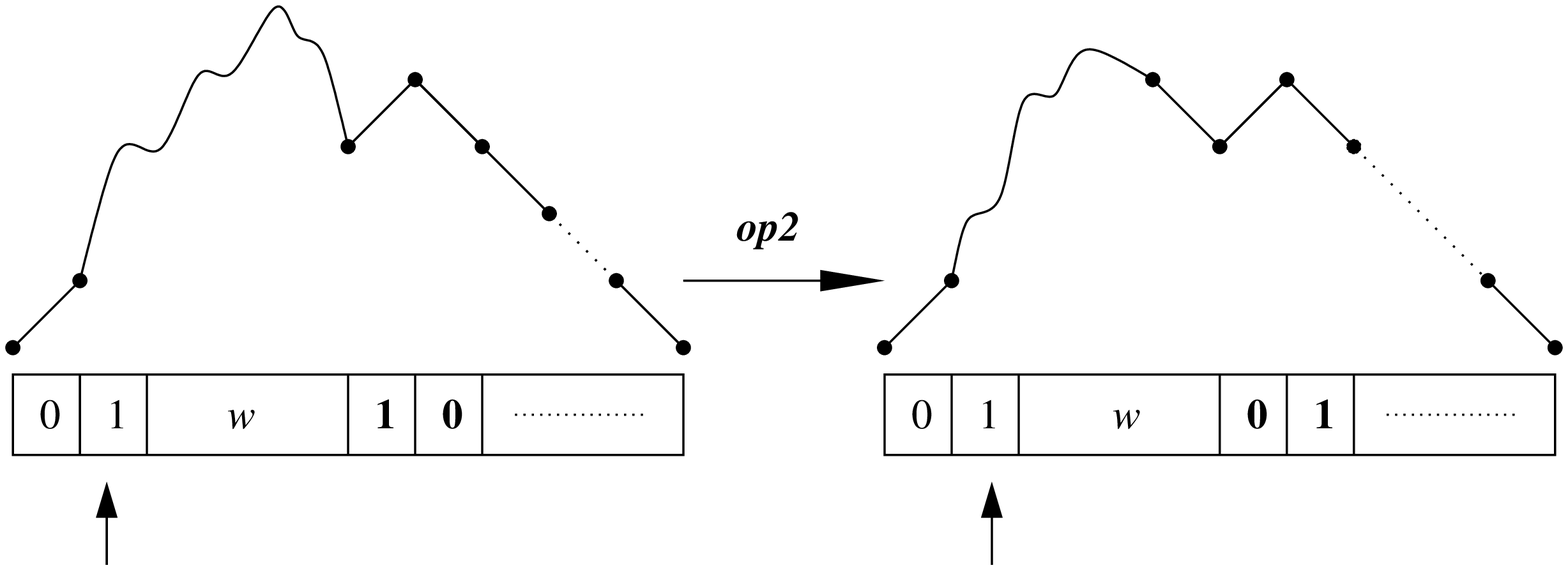}
\caption{Action of \emph{op2} operation on the array.}
\label{tr:op2v}
\end{centering}
\end{figure}

\begin{figure}[!ht]
\begin{centering}
\includegraphics[width=0.7\textwidth]{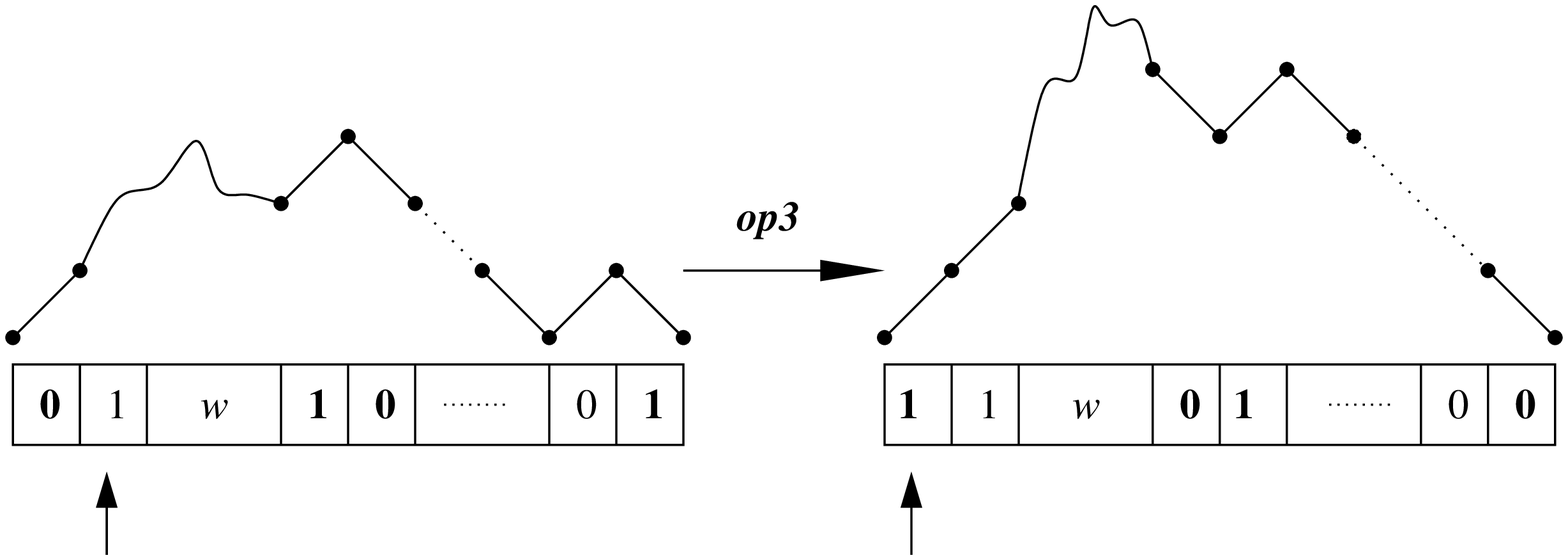}
\caption{How word of bits changes by \emph{op3}.}
\label{tr:op3v}
\end{centering}
\end{figure}

We pass to analyze \textbf{Algorithm 1}; our aim is to realize a
method which maintains constant the number of mean operations
while each object in $\mathcal{D}_n$ is generated. Now, if we
associate to each path a binary word by coding with 1 a northeast
step and with 0 a southeast, then it's clear that the three
operations are characterized by a constant number of actions which
exchange steps in the path. Indeed, we represent the word by a
circular array where the last position is followed by the first
one; we introduce a pointer to the first position of the array
which always corresponds to the first step of the path (see Figure
\ref{tr:vet}).
\newline Now, \emph{op1} is equivalent to exchange the
first bit 1 of the path with the first bit 0 of its last descent
and then to move forward the pointer one position (the action of
\emph{op1} on the array is illustrated in Figure \ref{tr:op1v}).
\newline \emph{Op2} is equivalent to exchange the bits of the last sequence $1 0$
in the array, while the pointer doesn't move (see Figure
\ref{tr:op2v}).
\newline Finally, \emph{op3} is
equivalent to exchange the bits of the last and second-last pairs
$10$ and then to move backward the pointer one position (see
Figure \ref{tr:op3v}).

It's clear that the three operations require a constant number of
actions independently of the length of the paths and
\textbf{Algorithm 1} is a constant amortized time (CAT) algorithm.

\section{Conclusions} We have presented a method to
generate all and only paths $\in \mathcal{D}_n$. The practical
advantages of our method are that it uses directly the
combinatorial objects and it generates all the paths $\in
\mathcal{D}_n$, with fixed $n$, without using the objects with
smaller size. So, as we have seen, our work presents
two different strategies 
which are closely connected. Indeed, the former can be described
by a rooted tree and the latter uses three operations, for listing
the objects, which are equivalent to visit this tree. Moreover, we
have proved that \textbf{Algorithm 1} is a constant amortized time
algorithm since it uses only a constant number of computations per
object.

Our studies have proved that the basic idea of this algorithm
allows to obtain similar results for other classes of paths like
Grand Dyck ($\mathcal{G}_n$) and Motzkin ($\mathcal{M}_n$) paths;
indeed, it's possible to obtain all the paths of $\mathcal{G}_n$
or $\mathcal{M}_n$ using operations very similar to \emph{op1},
\emph{op2} and \emph{op3}.


Moreover, it is reasonable to think that this method could be
applicable to other kinds of paths or to other combinatorial
classes which are in bijection with the studied paths. For example
we could study the classes of polyominoes or permutations
enumerated by Catalan, Motzkin or Gran Dyck numbers (for
definitions see for example \cite{21}).

\end{document}